\documentclass[11pt]{amsart}

\newtheorem{theorem}{Theorem}[section]
\newtheorem{proposition}{Proposition}[section]
\newtheorem{lemma}{Lemma}[section]

\theoremstyle{definition}
\newtheorem{remark}{Remark}[section]

\numberwithin{equation}{section}

\begin{document}

\title[Criterion for Riemann Hypothesis] {A new necessary and sufficient condition for
the Riemann hypothesis}

\author{Luis B\'{a}ez-Duarte}

\date{13 July 2003}
\email{lbaezd@cantv.net}


\begin{abstract}
We give a new equivalent condition for the Riemann hypothesis consisting in an
order condition for certain finite rational combinations of the values of
$\zeta(s)$ at even positive integers.
\end{abstract}

\maketitle

\section{Introduction and Preliminairies}
In this note we shall prove the following theorem:

\begin{theorem}\label{equivthm}
Let

\begin{equation}\label{coeffs}
c_k:=\sum_{j=0}^k (-1)^j {k \choose j} \frac{1}{\zeta(2j+2)},
\end{equation} 
then the Riemann hypothesis is true if and only if

\begin{equation}\label{condition}
c_k\ll k^{-\frac{3}{4}+\epsilon}, \ \ \ (\forall \epsilon>0).
\end{equation}
\end{theorem}

\begin{remark}
It will be seen below that unconditionally

\begin{equation}
c_k\ll k^{-\frac{1}{2}}.
\end{equation}
\end{remark}

\begin{remark}
It is quite obvious how one can trivially modify the proof of the theorem to obtain a
more general result:

\begin{theorem}\label{equivthm2}
A necessary and sufficient condition for
$\zeta(s)\not=0$ in the half-plane $\Re(s)>2(1-\alpha)$ is

\begin{equation}\label{condition}
c_k\ll k^{-\alpha +\epsilon}, \ \ \ (\forall \epsilon>0).
\end{equation}
\end{theorem}
\noindent However we shall eschew such gratuitous generalizing at this stage.
 \end{remark}
\ \\
Necessary and sufficient conditions for the Riemann hypothesis depending only on
values of $\zeta(s)$ at positive integers have been known for a long time, e.g.
those of M. Riesz \cite{riesz} and Hardy-Littlewood \cite{hardy}. M. Riesz's
criterion, for example, states that the Riemann hypothesis is true if and only if

$$
\sum_{k=1}^\infty\frac{(-1)^{k+1}x^k}{(k-1)!\zeta(2k)}
=O(x^{\frac{1}{4}+\epsilon}), \ \ \ (x\rightarrow +\infty).
$$
\ \\
We believe our condition is new and it is definitely simpler, as it only
involves finite rational combinations of the values $\zeta(2h)$, and seems well
posed for numerical calculations. This work however did not originate as an attempt to
simplify Riesz's criterion. It arose rather as a consequence of our note
\cite{baez} on Maslanka's expression of the Riemann zeta function
(\cite{maslanka1},\cite{maslanka2}) in the form

$$
(s-1)\zeta(s)=\sum_{k=0}^\infty A_k P_k\left(\frac{s}{2}\right).
$$
Here the $P_k(s)$ are the \textit{Pochhammer polynomials}

\begin{equation}\label{poch}
P_k(s):=\prod_{r=1}^k \left(1-\frac{s}{r}\right),
\end{equation}
which will appear prominently in the proof of Theorem \ref{equivthm}. Two 
elementary facts about them shall be needed: firstly

\begin{equation}\label{pochcomb}
(-1)^k{\frac{s}{2}-1 \choose k}=P_k\left(\frac{s}{2}\right),
\end{equation}
which is essentially a matter of notation, and secondly a standard estimate
given here without proof:

\begin{lemma}\label{pochestlemma} 
For every circle $|s|<r<\infty$ there is a positive constant $C_r$ such that

\begin{equation}\label{pochest}
|P_k(s)|\leq C_r k^{-\Re(s)}.
\end{equation}
\end{lemma}

\section{Sufficiency of the condition}
The sufficiency of the condition (\ref{condition})  follows from writing
$(\zeta(s))^{-1}$ as a series of Pochhammer polynomials.

\begin{proposition}[Sufficiency of the condition]\label{zetainv}
If $c_k\ll k^{-\frac{3}{4}+\frac{1}{2}\epsilon}$ for any $\epsilon>0$, then

\begin{equation}\label{zetainvrepr}
\frac{1}{\zeta(s)}=\sum_{k=0}^\infty c_k P_k\left(\frac{s}{2}\right), \
\ \ \ (\Re(s) > \frac{1}{2}),
\end{equation}
where the series converges uniformly on compact subsets of the half-plane.
\end{proposition}

\begin{remark}
Since it shall be shown that actually $c_k\ll k^{-\frac{1}{2}}$ it follows
modifying trivially the above argument that the representatiom (\ref{zetainvrepr}) for
$(\zeta(s))^{-1}$ is \textit{unconditionally} valid at least in the half-plane
$\Re(s)>1$. 
\end{remark}

We need a lemma before proving Proposition \ref{zetainv}.

\begin{lemma}\label{zetacoeffs}
Define

\begin{equation}\label{zetacoeffs1}
q_k:=\sum_{n=1}^\infty \frac{1}{n^2}\left(1-\frac{1}{n^2}\right)^k,
\end{equation}
then

\begin{equation}\label{zetacoeffsest}
q_k \ll k^{-\frac{1}{2}}.
\end{equation}
\end{lemma}

\begin{proof}
Let $\overline{B_1}(x)=x-[x]-\frac{1}{2}$. By the Euler-MacLaurin formula we
have for
$k\geq 1$

\begin{eqnarray}\label{qk}
q_k
&=&
\int_{1}^\infty \frac{1}{x^2} \left(1-\frac{1}{x^2}\right)^k dx
+\int_{1}^\infty \overline{B_1}(x)\frac{d}{dx}
\frac{1}{x^2} \left(1-\frac{1}{x^2}\right)^k dx\\\nonumber
&=&
\sqrt{\pi}\frac{\Gamma(k+1)}{\Gamma(k+\frac{3}{2})}
+\int_{1}^\infty \overline{B_1}(x)\frac{d}{dx}\left(
\frac{1}{x^2} \left(1-\frac{1}{x^2}\right)^k\right) dx.
\end{eqnarray}
Clearly

\begin{equation}\label{int1}
\frac{\Gamma(k+1)}{\Gamma(k+\frac{3}{2})}\ll k^{-\frac{1}{2}},
\end{equation}
and, letting $V(f(x))$ denote the total variation
 of $f(x)$ in $[1,\infty)$, we see that

\begin{eqnarray}\label{rem1}\nonumber
\left|\int_{1}^\infty \overline{B_1}(x)\frac{d}{dx}
\left(\frac{1}{x^2} \left(1-\frac{1}{x^2}\right)^k\right) dx\right|
&\ll&
\int_{1}^\infty 
\left|\frac{d}{dx} 
\left(\frac{1}{x^2} \left(1-\frac{1}{x^2}\right)^k\right)\right|dx
\\\nonumber
&=&
V\left(\frac{1}{x^2}\left(1-\frac{1}{x^2}\right)^k\right)\\\nonumber
&=&
2 \max_{1\leq x < \infty} \frac{1}{x^2}\left(1-\frac{1}{x^2}\right)^k
\\
&=&
\frac{2}{k+1}\left(1-\frac{1}{k+1}\right)^k \ll k^{-1}.
\end{eqnarray}
Hence, (\ref{qk}),(\ref{int1}) and (\ref{rem1}) achieve (\ref{zetacoeffsest}).
\end{proof}

\begin{proof}[Proof of Proposition \ref{zetainvrepr}]
First note that

\begin{eqnarray}\label{coeffs1}
c_k
&=&\
\sum_{j=0}^k (-1)^j {k \choose j} \frac{1}{\zeta(2j+2)}\\\nonumber
&=&
\sum_{j=0}^k (-1)^j {k \choose j} \sum_{n=1}^\infty \frac{\mu(n)}{n^{2j+2}}
\\\nonumber
&=&
\sum_{n=1}^\infty\frac{\mu(n)}{n^2}
\sum_{j=0}^k (-1)^j {k \choose j}\frac{1}{n^{2j}}\\\nonumber
&=&
\sum_{n=1}^\infty\frac{\mu(n)}{n^2}\left(1-\frac{1}{n^2}\right)^k.
\end{eqnarray}

Starting now with $\Re(s)>1$ we have

\begin{eqnarray}\label{eqx}\nonumber
\frac{1}{\zeta(s)}=\sum_{n=1}\frac{\mu(n)}{n^s}
&=&
\sum_{n=1}^\infty\frac{\mu(n)}{n^2}\left(\frac{1}{n^2}\right)^{\frac{s}{2}-1}
\\\nonumber
&=&
\sum_{n=1}^\infty\frac{\mu(n)}{n^2}
\left(1-\left(1-\frac{1}{n^2}\right)\right)^{\frac{s}{2}-1}
\\\nonumber
&=&
\sum_{n=1}^\infty\frac{\mu(n)}{n^2}
\sum_{k=0}^\infty (-1)^k {\frac{s}{2}-1 \choose k} \left(1-\frac{1}{n^2}\right)^k
\\
&=&
\sum_{n=1}^\infty\frac{\mu(n)}{n^2}
\sum_{k=0}^\infty P_k\left(\frac{s}{2}\right)\left(1-\frac{1}{n^2}\right)^k.
\end{eqnarray}
These summations can be interchanged because calling

$$
S=\sum_{n=1}^\infty\sum_{k=0}^\infty\frac{1}{n^2}
 \left|P_k\left(\frac{s}{2}\right)\right|\left(1-\frac{1}{n^2}\right)^k,
$$
we see from Lemma \ref{pochestlemma} and Lemma \ref{zetacoeffs1} that

\begin{eqnarray}\nonumber
S 
&=&
\sum_{k=0}^\infty \left|P_k\left(\frac{s}{2}\right)\right|
\sum_{n=1}^\infty\frac{1}{n^2}
\left(1-\frac{1}{n^2}\right)^k\\\nonumber
&=&
\sum_{k=0}^\infty \left|P_k\left(\frac{s}{2}\right)\right|q_k
\ll \sum_{k=1}^\infty k^{-\frac{\Re(s)}{2} -\frac{1}{2}} < \infty.
\end{eqnarray}
Thus we proceed to interchange summations in (\ref{eqx}), taking into account
(\ref{coeffs1}), to obtain unconditionally for $\Re(s)>1$,

\begin{eqnarray}\label{uncond}\nonumber
\frac{1}{\zeta(s)}
&=&
\sum_{k=0}^\infty P_k\left(\frac{s}{2}\right)
\sum_{n=1}^\infty \frac{\mu(n)}{n^2}\left(1-\frac{1}{n^2}\right)^k\\
&=&
\sum_{k=0}^\infty c_k P_k\left(\frac{s}{2}\right).
\end{eqnarray}
But Lemma \ref{pochestlemma} together with the hypothesis
$c_k<<k^{\frac{3}{4}+\frac{1}{2}\epsilon}$ implies that the above series converges
uniformly on compacts of the half-plane $\Re(s)>\frac{1}{2}+\epsilon$. This means that
the series extends $(\zeta(s))^{-1}$ analytically to the half-plane
$\Re(s)>\frac{1}{2}$. 
\end{proof}

\section{Necessity of the condition}
\begin{proof}[Proof of the necessity of the condition]
Assume now that the Riemann hypothesis is true. If as usual we write

$$
M(x):=\sum_{n\leq x} \mu(n),
$$
we then have

$$
M(x)\ll x^{\frac{1}{2}+2\epsilon}, \ \ \ \ (\forall \epsilon>0).
$$
We can transform the second expression for $c_k$ in (\ref{coeffs1}) summing it
by parts to obtain

\begin{eqnarray}\nonumber
c_k
&=&
\int_1^\infty M(x)
\frac{d}{dx}\left(\frac{1}{x^2}\left(1-\frac{1}{x^2}\right)^k\right) dx
\\\nonumber
&=&
2\int_0^1 M\left(\frac{1}{x}\right)(1-x^2)^k((k+2)x^3 - x)dx.
\end{eqnarray}
Therefore

$$
|c_k|\leq 2(k+2)\int_0^1\left|M\left(\frac{1}{x}\right)\right|x^3(1-x^2)^kdx +
\int_0^1 \left|M\left(\frac{1}{x}\right)\right|x(1-x^2)^kdx,
$$
but (on the Riemann hypothesis)

$$
M\left(\frac{1}{x}\right)<<x^{-\frac{1}{2}-2\epsilon}, \ \ \ \ (x\downarrow 0),
$$
so that

\begin{equation}\label{ckest}
c_k \ll k\int_0^1 x^{\frac{5}{2}-2\epsilon}(1-x^2)^k dx +
\int_0^1 x^{\frac{1}{2}-2\epsilon}(1-x^2)^k dx.
\end{equation}
On the other hand, for $\Re(\lambda)>-1$ a classical beta integral result is

$$
\int_0^1 x^\lambda(1-x^2)^k dx =
\Gamma\left(\frac{\lambda+1}{2}\right)
\frac{\Gamma(k+1)}{\Gamma(k+\frac{1}{2}(\lambda+3))}\ll
k^{-\frac{1}{2}-\frac{\lambda}{2}},
$$
so that (\ref{ckest}) becomes

$$
c_k \ll k^{-\frac{3}{4}+\epsilon}.
$$
\end{proof}

\section{Results of some calculations}
A test for the first $c_k$ up to $k=1000$ shows a very pleasant smooth curve
which, on the meager strength of so limited a calculation, would seem to
indicate that 

$$
c_k k^{\frac{3}{4}}\log^2 k
$$
tends to a finite limit in a very regular way.
\\
\ \\

\bibliographystyle{amsplain}
  
 \ \\
\ \\
\noindent Luis B\'{a}ez-Duarte\\
Departamento de Matem\'{a}ticas\\
Instituto Venezolano de Investigaciones Cient\'{\i}ficas\\
Apartado 21827, Caracas 1020-A\\
Venezuela\\
\email{lbaezd@cantv.net}

\end{document}